\documentclass[10pt, a4paper]{amsart}
\usepackage[top=1.25in, bottom=1.25in, left=1.25in, right=1.25in]{geometry}
\usepackage{graphicx}
\usepackage{amsfonts}
\usepackage{epsf}
\usepackage{amssymb}
\usepackage{amsmath}
\usepackage{amscd}
\usepackage{amsthm}
\usepackage{tikz}
\usetikzlibrary{
  cd,
  calc,
  positioning,
  arrows,
  decorations.pathreplacing,
  decorations.markings,
}
\usepackage{pdfpages}
\usepackage{fancyhdr}
\usepackage{setspace}
\usepackage{hyperref}
\usepackage[all]{xy}
\usetikzlibrary{matrix}
\usepackage{verbatim}
\usepackage{enumerate}
\usepackage{mathrsfs}
\usepackage{pinlabel}
\usepackage{lscape}
\usepackage{color}
\usepackage{subfigure}
\usepackage[normalem]{ulem}

\newtheorem{theorem}{Theorem}[section]

\newtheorem{proposition}[theorem]{Proposition}
\newtheorem{corollary}[theorem]{Corollary}
\newtheorem{lemma}[theorem]{Lemma}

\newtheorem*{rep@theorem}{\rep@title}
\newcommand{\newreptheorem}[2]{
\newenvironment{rep#1}[1]{
\def\rep@title{#2 \ref{##1}}
\begin{rep@theorem}}
{\end{rep@theorem}}}
\newreptheorem{theorem}{Theorem}
\newreptheorem{lemma}{Lemma}
\newreptheorem{proposition}{Proposition}
\newreptheorem{corollary}{Corollary}

\newtheorem{innercustomgeneric}{\customgenericname}
\providecommand{\customgenericname}{}
\newcommand{\newcustomprop}[2]{%
  \newenvironment{#1}[1]
  {%
   \renewcommand\customgenericname{#2}%
   \renewcommand\theinnercustomgeneric{##1}%
   \innercustomgeneric
  }
  {\endinnercustomgeneric}
}
\newcustomprop{customprop}{Proposition}

\theoremstyle{definition}
\newtheorem{definition}[theorem]{Definition}
\newtheorem{example}[theorem]{Example}
\newtheorem{remark}[theorem]{Remark}
\newtheorem{problem}[theorem]{Problem}

\newtheorem*{case2'}{Case 2$'$}
\newtheorem{theorem-named}{}
\newtheorem{theorem-labeled}{Theorem}
\newtheorem{definition-named}{}
\newtheorem{conjecture-named}{}
\newtheorem{case-named}{}

\numberwithin{equation}{section}

\newcommand{\Z}{\mathbb{Z}}
\newcommand{\Q}{\mathbb{Q}}

\newcommand{\s}{\mathfrak{s}}

\def\Z{\mathbb{Z}}

\def\Q{\mathbb{Q}}

\def\hat{\widehat}

\DeclareMathOperator\Ker{Ker}

\makeatletter
\def\gt{g^\mathrm{top}_4}
\def\gz{g_{\Z}}
\def\gc{g_{\Z}^\mathrm{c}}
\def\gtc{g_{\Z}^\mathrm{t.c.}}
\newcommand{\Imm}{\text{Im}}

\begin{document}
\title{A note on the concordance $\Z$-genus}

\author{Allison N.\ Miller}
\thanks{The first author was partially supported by NSF grant DMS-1902880.}
\address{Department of Mathematics, Rice University, Houston, TX, United States}
\email{allison.miller@rice.edu}

\author{JungHwan Park}
\address{Department of Mathematical Sciences, KAIST, Daejeon, South Korea}
\email{jungpark0817@kaist.ac.kr}

\def\subjclassname{\textup{2020} Mathematics Subject Classification}
\expandafter\let\csname subjclassname@1991\endcsname=\subjclassname
\expandafter\let\csname subjclassname@2000\endcsname=\subjclassname

\begin{abstract}
We show that the difference between the topological 4-genus of a knot and the minimal genus of a surface bounded by that knot that can be decomposed into a smooth concordance followed by an algebraically simple locally flat surface can be arbitrarily large. This  extends work of Hedden-Livingston-Ruberman showing that there are 
topologically slice knots which are not smoothly concordant to any  knot with trivial Alexander polynomial.

\end{abstract}
\maketitle

\section{Introduction}\label{sec:intro}
The \emph{$\Z$-genus} of a knot $K$ in $S^3$, denoted by $\gz(K)$, is the minimal genus of an oriented properly embedded locally flat surface $\Sigma$ in the 4-ball $B^4$ such that the boundary of $\Sigma$ is $K$ and $\pi_1 (B^4 \smallsetminus \Sigma) \cong \mathbb{Z}$. Work of Freedman famously implies that a knot $K$ has $\gz(K)=0$ if and only if $K$ has trivial Alexander polynomial~\cite{Freedman:1982-1, Freedman-Quinn:1990-1, Garoufalidis-Teichner:2004-1}; later work of Feller  generalizes Freedman's result to show that $\gz(K)$ is bounded above by half the degree of the Alexander polynomial of $K$~\cite{Feller:2016-1}. The $\Z$-genus of $K$ can therefore be thought of as an algebraically controlled upper bound on the topological 4-ball genus of $K$, i.e. the minimal genus of \textit{any} oriented properly embedded locally flat surface in the 4-ball with boundary $K$ ; see Feller-Lewark for a precise statement of this fact~\cite{Feller-Lewark:2018-1,Feller-Lewark:2019-1}. Further, we define the \emph{$($smooth$)$ concordance $\Z$-genus} of a knot $K$, denoted by $\gc(K)$, to be the minimum value of $\gz(J)$ among all knots smoothly concordant to $J$. That is, $$\gc(K):=\min\{\gz(J)\mid J\text{ is smoothly concordant to } K\} .$$ Observe that a knot $K$ is smoothly concordant to a knot with trivial Alexander polynomial if and only if $\gc(K)=0$, and that  we obtain the following inequalities for any knot $K$ immediately from the definitions:
$$\gt(K) \leq \gc(K) \leq \gz(K),$$
where $\gt(K)$ denotes the topological $4$-ball genus of $K$.

We remark that taking connected sums of a smoothly slice knot with nontrivial determinant produces knots with vanishing $\gc$ and arbitrarily large $\gz$ (see e.g.\ Proposition~\ref{prop:branchedlower}). In particular, this implies that the gap between $\gc$ and $\gz$ can be made arbitrarily large.
 In~\cite{Hedden-Livingston-Ruberman:2012-1}, Hedden-Livingston-Ruberman used the Heegaard Floer correction terms of 3-manifolds to show that there is a topologically slice knot $K$ which is not smoothly concordant to any knot with trivial Alexander polynomial, which in particular implies that $\gt(K) < \gc(K)$. We further show that the gap between these two invariants can be arbitrarily large.

\begin{theorem}\label{thm:main}There exist topologically slice knots with arbitrarily large concordance $\Z$-genus. More precisely, there exists a topologically slice knot $K_*$ such that for each $n \in \mathbb{N}$ its $n$-fold connected self-sum has concordance $\Z$-genus at least $n$. Furthermore, we can choose  $K_*$ so that it has smooth $4$-ball genus one.\end{theorem}


We outline the proof of the main result. We first construct a topologically slice knot $K_*$ with smooth 4-ball genus one in Example~\ref{exl:ourexl}. Let $n$ be a positive integer and $\Sigma_2(\#^n K_*)$ be the 2-fold branched cover of $S^3$ branched over $\#^n K_*$. We then observe that for any nontrivial 5-torsion element $x$ of $H_1(\Sigma_2(\#^n K_*))$, we may multiply an appropriate constant to obtain $z = c\cdot x$ so that the Heegaard Floer $\overline{d}$-invariant $\overline{d}(\Sigma_2(\#^nK_*), \s_z)$ is nonvanishing. We use this to show that there is no rational homology cobordism from $\Sigma_2(\#^n K_*)$ to any 3-manifold $Y$ with $H_1(Y)$ generated by at most $2n$ elements, and conclude that $\gc(\#^nK_*)\geq n$.

On a similar note, Hom~\cite{Hom:2015-1} showed that there are topologically slice knots with smooth 4-ball genus equal to one and arbitrarily large concordance genus. Recall that the concordance genus is greater than or equal to concordance $\mathbb{Z}$-genus since the fundamental group of the complement of a pushed in Seifert surface is equal to $\Z$.

Theorem~\ref{thm:main} is related  to the following difficult problem.
We can analogously define the \emph{topological concordance $\Z$-genus} of a knot $K$ as $$\gtc(K):=\min\{\gz(J)\mid J\text{ is topologically concordant to } K\} .$$ So for any knot $K$ we have that 
$$\gt(K) \leq \gtc(K) \leq \gc(K) \leq \gz(K).$$

\noindent Livingston~\cite{Livingston:2004-1} showed that the gap between $\gt$ and $\gtc$ can also be made arbitrarily large by using Casson-Gordon invariants (see Remark~\ref{rmk:livingston}), but the following remains open.

\begin{problem}\label{prob:smoothtop}
Show that topological $4$-ball genus can be strictly less than topological concordance smooth $4$-ball genus, i.e. 
find a knot $K$ with $\gt(K)=g \geq 1$ such that $K$ is not topologically concordant to any knot $J$ with $g_4(J) = g$, where $g_4(J)$ denotes the smooth $4$-ball genus of $J$.
\end{problem}

The following problem about properly embedded surfaces is a priori easier, but is nevertheless still open.
\begin{problem}\label{prob:smoothtopsurface}
Find an oriented properly embedded locally flat genus $g\geq 1$ surface $\Sigma$ in the 4-ball such that $\Sigma$ is not ambiently isotopic rel boundary to any surface $\Sigma'$ that can be decomposed into a topological concordance glued to a smoothly embedded surface.
\end{problem}

 Note that a positive solution to Problem~\ref{prob:smoothtop} certainly implies a positive solution to Problem~\ref{prob:smoothtopsurface}. 
While we avoid the technical details of topological 4-manifold theory needed to prove this formally, it is known both that every topologically slice knot is topologically concordant to the unknot, which of course has smooth $4$-ball genus equal to zero, and that every topologically embedded disk can be decomposed into a topological concordance glued to the standard slice disk for the unknot.  That is, the $g=0$ case of Problems~\ref{prob:smoothtop} and~\ref{prob:smoothtopsurface} is resolved.

\subsection*{Acknowledgments} 

This paper began in conversations during the first author's visit to the School of Mathematics at  Georgia Institute for Technology to speak in the Geometry Topology Seminar; she thanks Miriam Kuzbary and the School as a whole for the invitation and for their generous hospitality. We also thank Peter Feller, Mark Powell, and Charles Livingston for helpful conversations.

\subsection*{Notation and conventions}
In this paper, we work in the smooth category unless specified otherwise. Let $-K$ be the reverse of the mirror image of a knot $K$ and $-Y$ be the manifold $Y$ with reversed orientation. For any given knot $K$ in $S^3$,  the 2-fold branched cover of $S^3$ branched over $K$ is denoted by $\Sigma_2(K)$, 
 and for $n\in \Z$, the $3$-manifold obtained by $n$-framed Dehn surgery on $S^3$ along $K$ is denoted by $S^3_n(K)$.

\section{Proof of Theorem~\ref{thm:main}}\label{sec:pfthm}
Feller-Lewark defined the \emph{algebraic genus} of a knot $K$ in~\cite{Feller-Lewark:2018-1} and further showed that this invariant coincides with the $\Z$-genus~\cite{Feller-Lewark:2019-1}. This implies the following key  ingredient for the proof of Theorem~\ref{thm:main}. For a finite abelian group $G$, we denote by $r(G)$ the so-called ``generating rank'', i.e. the minimal number of generators of $G$.

\begin{proposition}[Proposition 12 of \cite{Feller-Lewark:2018-1}]
\label{prop:branchedlower} If $\Sigma_2(K)$ is the 2-fold branched cover of $S^3$ branched over a knot $K$, then
\[\pushQED{\qed}r(H_1(\Sigma_2(K)))\leq 2\gz(K).\qedhere\]
\end{proposition}

The other ingredient for the proof comes from an analysis of metabolizers for  torsion linking forms. Recall that there is a nonsingular symmetric linking form $$\lambda \colon H_1(Y) \times H_1(Y) \to \Q/\Z,$$ for a rational homology sphere $Y$. If $-Y$ is the 3-manifold obtained by reversing the orientation of $Y$, then the linking form on $H_1(-Y)$ is given by $-\lambda$. We say that a subgroup $M \subseteq H_1(Y)$ is a metabolizer if $M = M^{\perp}$ with respect to $\lambda$; note that for a metabolizer $M$ the nonsingularity of $\lambda$ immediately implies that $|M|^2=|H_1(Y)|$. 
It is well known that if $Y$ bounds a rational homology ball $W$, then the kernel of the inclusion-induced map $H_1(Y) \to H_1(W)$ is a metabolizer for $(H_1(Y),\lambda)$.  

For a finite abelian group $G$ and a prime $p$, let $G_p$ denote the $p$-primary subgroup of $G$, i.e.\ the subgroup consisting of all elements of $G$ which are $p^k$-torsion for some positive integer $k$.


\begin{lemma}\label{lem:key}
Let $Y_1$ and $Y_2$ be rational homology 3-spheres and let $p$ be a prime integer. Furthermore, for $i=1,2$, let $\lambda_i$ be the linking form on $H_1(Y_i)$, and let $M$ be a metabolizer for $(H_1(Y_1\#-Y_2), \lambda_1 \oplus -\lambda_2)$. 

If $H_1(Y_1)_p \cong(\Z/{p^2}\Z)^{2n}$ and $r(H_1(Y_2)_p) \leq 2m <2n$ for nonnegative integers $n$ and $m$, then $M \cap H_1(Y_1)$ contains a subgroup isomorphic to $(\Z/p\Z)^{n-m}$.\end{lemma}

\begin{proof}
First, note that $M$ splits as a direct sum of its $q$-primary components over all prime $q$, each of which is a metabolizer for the $q$-primary component of $H_1(Y_1 \# -Y_2)$. We can therefore assume without loss of generality that $H_1(Y_1)$ and $H_1(-Y_2)$ are equal to their $p$-primary components, and do so for ease of notation. 

For $i=1,2$, let $\pi_i \colon M \to H_1(Y_i)$ be the restriction of the projection map $$H_1(Y_1 \# -Y_2) \cong H_1(Y_1) \oplus H_1(-Y_2) \to H_1(Y_i).$$ Observe that $\Ker (\pi_2)= M \cap H_1(Y_1)$, which will be used later. We now argue that $p \Imm(\pi_2)$ has the property that $-\lambda_2|_{p \Imm(\pi_2) \times p \Imm(\pi_2)} = 0$. If $b_1, b_2 \in \Imm(\pi_2)$, then there exist $a_1, a_2 \in H_1(Y_1)$ such that $(a_1, b_1), (a_2, b_2) \in M$, and 
\[ (\lambda_1 \oplus - \lambda_2)( (a_1, b_1), (a_2, b_2))
=\lambda_1(a_1, a_2) - \lambda_2(b_1, b_2)=0 \in \Q/\Z.
\]
Since $H_1(Y_1)$ is annihilated by $p^2$, we have that
\[-\lambda_2(pb_1, pb_2)= -p^2 \lambda_2(b_1, b_2)= -p^2 \lambda_1(a_1, a_2) =  -\lambda_1(p^2a_1, a_2) = 0 \in \Q/\Z,\] as desired. Moreover, since the linking form $-\lambda_2$ is nonsingular, we have 
\begin{equation}\label{eq:1}
|p \Imm(\pi_2)| \leq |H_1(-Y_2)|^{1/2}.
\end{equation}

Recall that since we are assuming $r(H_1(Y_2)) \leq 2m$ and since $\Imm(\pi_2)$ is a subgroup of $H_1(-Y_2)$, we have that
\[\Imm(\pi_2) \cong \bigoplus_{k=1}^{\ell}  (\Z/ p^k\Z)^{n_k} 
\quad\text{and}\quad
p \Imm(\pi_2) \cong \bigoplus_{k=1}^{\ell}  ( \Z/ p^{k-1}\Z)^{n_k} ,
\quad\text{where}\quad  \sum_{k=1}^{\ell} n_k \leq 2m.
\]
Therefore, 
\[ \frac{|\Imm(\pi_2)|}{|p \Imm(\pi_2)|}= \frac{ p^{\sum_{k=1}^{\ell} k n_k}}{p^{\sum_{k=1}^{\ell} (k-1)n_k}}= p^{\sum_{k=1}^{\ell} n_k}\ \leq p^{2m}.
\]
So, by combining inequality \eqref{eq:1} with the above inequality, we conclude that $$|\Imm(\pi_2)| \leq p^{2m} |H_1(-Y_2)|^{1/2}.$$
Also, recall that $|M| = |H_1(Y_1)|^{1/2} |H_1(-Y_2)|^{1/2}$ and $\Ker(\pi_2)= M \cap H_1(Y_1)$. Therefore,
\[ |M \cap H_1(Y_1)|= |\Ker(\pi_2)| = \frac{|M|}{| \Imm(\pi_2)|} \geq 
\frac{|H_1(Y_1)|^{1/2} |H_1(-Y_2)|^{1/2}}{p^{2m} |H_1(-Y_2)|^{1/2}}= p^{2n-2m}.
\]
Finally, since $M \cap H_1(Y_1)$ is a subgroup of $(\Z/p^2\Z)^{2n}$, we conclude that $$M \cap H_1(Y_1) \cong (\Z/p\Z)^{k_1} \oplus (\Z/p^2\Z)^{k_2}, \quad\text{where}\quad k_1 + 2 k_2 \geq 2n-2m.$$ 
Therefore,
\[ 2n-2m \leq k_1 + 2 k_2 \leq 2k_1 + 2k_2,\] and we conclude that $M \cap H_1(Y_1)$ contains a subgroup isomorphic to $(\Z/p\Z)^{n-m}$.\end{proof}


Now, we collect some necessary background on the Heegaard-Floer correction term $d(Y,\mathfrak{s})\in \mathbb{Q}$ associated to a rational homology sphere $Y$ with a Spin$^c$ structure~$\mathfrak{s}$~\cite{Ozsvath-Szabo:2003-2}. We first recall the following definition from~\cite{Hedden-Livingston-Ruberman:2012-1}.

\begin{definition}[{\cite[Definition 2.2 and 3.1]{Hedden-Livingston-Ruberman:2012-1}}]Let $Y$ be a $\Z/2\Z$-homology sphere. For $z\in H_1(Y)$, let $\mathfrak{s}_z$ be the unique Spin$^c$ structure of $Y$ which satisfies $c_1(\mathfrak{s}_z)=2\hat{z}\in H^2(Y)$ where $\hat{z}$ is the Poincar\'{e} dual of $z$. In particular, $\s_0$ is the unique Spin structure on $Y$. Lastly, we define $\overline{d}(Y,\mathfrak{s}_z):=d(Y,\s_z)-d(Y,\s_0).$\end{definition}

Further, the correction term is additive under connected sums \cite[Theorem 4.3]{Ozsvath-Szabo:2003-1} and is a Spin$^c$ rational homology cobordism invariant \cite[Theorem 1.2]{Ozsvath-Szabo:2003-1}. Hence we get the following lemma (see e.g.\ {\cite[Proposition 2.1]{Hedden-Livingston-Ruberman:2012-1}}).

\begin{lemma}\label{lem:dinvt}If $K$ and $J$ are two concordant knots, then there exists a metabolizer $M \subseteq H_1(\Sigma_2(K)\# - \Sigma_2(J))$ such that for each element $(m_K, m_J) \in M$ where $m_K \in H_1(\Sigma_2(K))$ and $m_J \in H_1(- \Sigma_2(J))$, we have that $d(\Sigma_2(K),\s_{m_K})+ d(-\Sigma_2(J),\s_{m_J}) =0.$\qed
\end{lemma}



We are now ready to prove the following corollary.

\begin{corollary}\label{cor:mainobstr}
Let $K$ be a knot with $\gc(K) = m$ and $p$ be a prime integer. If $H_1(\Sigma_2(K))_p \cong (\Z/p^2\Z)^{2n}$ and $n> m$, then there exists a subgroup $N \subseteq H_1(\Sigma_2(K))_p$ such that $N \cong (\Z/p\Z)^{n-m}$ and $\overline{d}(\Sigma_2(K), \s_z)=0$ for each $z \in N$. 
\end{corollary}

\begin{proof} Suppose $K$ is concordant to a knot $J$ with $\gz(J) = m$. Then by Proposition~\ref{prop:branchedlower} $r(H_1(\Sigma_2(J))) \leq 2m$. As in the proof of Lemma~\ref{lem:key} we assume that $H_1(\Sigma_2(K))$ and $H_1(-\Sigma_2(J))$ are equal to their $p$-primary components. 

By Lemma~\ref{lem:dinvt}, there exists a metabolizer $M \subseteq H_1(\Sigma_2(K)\# - \Sigma_2(J))$ such that for each element $(m_K, m_J) \in M$ where $m_K \in H_1(\Sigma_2(K))$ and $m_J \in H_1(- \Sigma_2(J))$,
\begin{equation}\label{eq:2}
d(\Sigma_2(K),\s_{m_K})+ d(-\Sigma_2(J),\s_{m_J}) =0.\end{equation}

Lastly, since $\Sigma_2(K)$ and $\Sigma_2(J)$ satisfy the assumption of Lemma~\ref{lem:key}, we may conclude that $M \cap H_1(Y_1)$ contains a subgroup isomorphic to $(\Z/p\Z)^{n-m}$. Set $N:=M \cap H_1(Y_1)$, then note that for each $z \in N$, we have $(0,0), (z,0) \in N$. The proof is complete by using equation~\eqref{eq:2}.
\end{proof}


\begin{example}\label{exl:ourexl}
Our examples come from the cabling construction. Let $K_{p,q}$ denote the $(p, q)$-cable of $K$, where $p$ is the longitudinal winding. We define our knot to be the following
$$K_*:= D_{2,25} \# -D_{2,23} \# -T_{2,25}\#  T_{2,23},$$
where $D$ is the $4$-fold connected self-sum of the positive Whitehead double of the right-handed trefoil knot. Since $D$ is topologically slice, $D_{p,q}$ is topologically concordant to the torus knot $T_{p,q}$. Therefore $K_*$ is topologically concordant to 
$$K_*':=T_{2,25} \# -T_{2,23} \# -T_{2,25}\#  T_{2,23},$$
hence is topologically slice. 
Also, note that by performing two crossing changes with opposite signs on $K_*$ (one in the $ D_{2,25}$ summand and one in the $-T_{2, 25}$ summand), 
we obtain the smoothly slice knot
$D_{2,23} \# -D_{2,23} \# -T_{2,23}\#  T_{2,23}.$ Hence $K_*$ bounds a smoothly embedded genus one surface in the 4-ball. Finally, by \cite{Akbulut-Kirby:branched-covers}, the 2-fold branched cover of $S^3$ branched over $K_*$ is given by
$$\Sigma_2(K_*)\cong S^3_{25}(D \# D^r) \#  -S^3_{23}(D \# D^r)\# -S^3_{25}(U)\# S^3_{23}(U),$$ where $D^r$ is a knot obtained by reversing the orientation of $D$. \end{example}

Next, we compute the correction terms for the 2-fold branched cover. Recall that for any knot $K$ in $S^3$, there are nonnegative integer-valued smooth concordnace invariants $V_i(K)$, introduced by Rasmussen~\cite{Rasmussen:2003-1}. Furthermore, Ni and Wu~\cite{Ni-Wu:2015-1} showed that the correction terms of all Dehn surgeries on $K$ can be computed from these invariants.

\begin{proposition}[{\cite[Proposition 1.6 and Remark 2.10]{Ni-Wu:2015-1}}]\label{proposition:NiWu} 
Let $K$ be a knot and $U$ be the unknot. If $n$ is a positive integer, then
\[\pushQED{\qed}d(S^3_n(K),\s_{i}) = d(S^3_n(U),\s_{i})-2\max \{V_{i}(K),V_{n-i}(K)\}.\qedhere\]
\end{proposition}

Furthermore, the correction terms for Dehn surgeries on the unknot are computed in \cite{Ozsvath-Szabo:2003-1}.

\begin{proposition}[{\cite[Proposition 4.8]{Ozsvath-Szabo:2003-1}}]\label{prop:dlens} If $U$ is the unknot and $n$ is a positive integer, then
\[\pushQED{\qed}d(S^3_n(U),\s_{i}) = \frac{(n-2i)^2}{4n}-\frac{1}{4}.\qedhere\]
\end{proposition}

Hom and Wu~\cite{Hom-Wu:2016-1} introduced a nonnegative integer valued smooth concordance invariant~$\nu^+$. Moreover, for any knot $K$, we have that $\nu^+(K)=0$ if and only if $V_0(K)=0$ \cite[Proposition 2.3]{Hom-Wu:2016-1}. Following~\cite{Kim-Park:2018-1} (see also \cite{Hom:2017-1}), we say that two knots $K$ and $K'$ are \emph{$\nu^+$-equivalent} if $$\nu^+(K \# - K') =\nu^+(K' \# - K)=0,$$ and it forms a equivalence relation on the set of concordance classes of knots. It is well-known that two $\nu^+$-equivalent knots have the same $V_i$ invariants (see e.g.\ \cite[Proposition 3.11]{Kim-Kracatovich-Park:2019-1}). We will use the fact that $D \# D^r$ is $\nu^+$ equivalent to the torus knot $T_{2,17}$. This fact might be well-known to experts, but we sketch the proof here for the reader’s convenience. 

\begin{proposition}\label{prop:DDT217}The knot  $D \# D^r$ in Example~\ref{exl:ourexl} is $\nu^+$-equivalent to the torus knot $T_{2,17}$. In particular, we have
$$V_i(D \# D^r)= \left\{
\begin{array}{cl}4&\text{for }i=0,\\
2&\text{for }i=5,\\
0&\text{for }i\geq 8.
\end{array}\right.$$
\end{proposition}
\begin{proof} If $K_i$ and $K'_i$ are $\nu^+$-equivalent for $i = 0, 1$, then $K_0\#K_1$ and $K'_0\#K'_1$ are also $\nu^+$-equivalent (see e.g.\ \cite[Proposition 3.12]{Kim-Kracatovich-Park:2019-1}). Also, $D^r$ is $\nu^+$-equivalent to $D$ since their knot Floer complexes are chain homotopy equivalent \cite[Proposition 3.9]{ozsvath2004holomorphicknot}. Hence we have that $D \# D^r$ is $\nu^+$-equivalent to the $8$-fold connected self-sum of the positive Whitehead double of the right-handed trefoil knot. Lastly, we apply a result of Hedden, Kim, and Livingston~\cite[Proposition 6.1]{Hedden-Kim-Livingston:2016-1} which states that the $n$-fold connected self-sum of the positive Whitehead double of the right-handed trefoil knot is $\nu^+$-equivalent to $T_{2,2n+1}$. This proves the first part of the statement.

The last part follows from two facts. First, as mentioned above, two $\nu^+$-equivalent knots have the same $V_i$ invariants. The second fact is that $V_i$ invariants for any alternating knot $K$ are determined by their signatures as follows (see e.g.\ \cite[Theorem 2]{Hom-Wu:2016-1} and \cite[Section 5.1]{Aceto-Golla:2017-1}):
$$V_i(K)=\max\left\{\left \lceil{-\frac{\sigma(K)+2i}{4}}\right \rceil,0  \right\} .$$
This completes the proof.\end{proof}

Now, we are ready to compute the correction terms.
\begin{proposition}\label{prop:d-of-2fold}Let $K_*$ be the knot in Example~\ref{exl:ourexl} and $\#^nK_*$ be the $n$-fold connected self-sum of $K_*$. If $N\subseteq H_1(\Sigma_2(\#^nK_*))_5$ is a nontrivial subgroup, then there exists an element $z \in N$ such that $\overline{d}(\Sigma_2(\#^nK_*), \s_z) \neq 0$.\end{proposition}
\begin{proof} By Example~\ref{exl:ourexl}, we have $\Sigma_2(K_*)\cong S^3_{25}(D \# D^r) \#  -S^3_{23}(D \# D^r)\# -S^3_{25}(U)\# S^3_{23}(U)$ and $H_1(\Sigma_2(K_*))_5\cong \Z/25\Z \oplus \Z/25\Z \cong H_1(S^3_{25}(D \# D^r) \# -S^3_{25}(U))$. Note that by the additivity of the correction terms, we have that 
\[
\overline{d}(\Sigma_2(K_*), \s_{(i,j)})= \overline{d}(\Sigma_2(K_*), \s_{(i,0)})+\overline{d}(\Sigma_2(K_*), \s_{(0,j)})
.\]

For $(i,0) \in H_1(\Sigma_2(K_*))_5$, we have by Proposition~\ref{proposition:NiWu}, Proposition~\ref{prop:dlens}, and Proposition~\ref{prop:DDT217} that
\begin{align*}
\overline{d}(\Sigma_2(K_*), \s_{(i,0)}) & = \overline{d}(S^3_{25}(U),\s_{i})-2\max \{V_{i}(D \# D^r),V_{25-i}(D \# D^r)\} +2V_0(D \# D^r) \\
  &= \left\{
\begin{array}{cl}0&\text{for }i=5 \text{ and } 20,\\
2&\text{for }i=10 \text{ and }15.
\end{array}\right.
\end{align*}
Similarly, for $(0,j) \in H_1(\Sigma_2(K_*))_5$ we have that
$$\overline{d}(\Sigma_2(K_*), \s_{(0,j)})=-\overline{d}(S^3_{25}(U),\s_{j})= \left\{
\begin{array}{cl}4&\text{for }j=5 \text{ and } 20,\\
6&\text{for }j=10 \text{ and }15.
\end{array}\right.$$ 

By additivity, we have that
 $\overline{d}(\Sigma_2(K_*), \s_x) \geq 0$ for any order $5$ element $x\in H_1(\Sigma_2(K_*))_5$. Moreover, at least one of $\overline{d}(\Sigma_2(K_*), \s_x)$ and $\overline{d}(\Sigma_2(K_*), \s_{2x})$ must be strictly positive.
 
 Now, let $N\subseteq H_1(\Sigma_2(\#^nK_*))_5$ be a nontrivial subgroup, and let $y$ be an order 5 element in $N$. Note that since $\Sigma_2(\#^n K_*) \cong  \#^n \Sigma_2(K_*)$, we can naturally write $y=(y_1,\dots, y_n)$ for $y_i \in H_1(\Sigma_2(K_*))_5$ and have that
 \[\overline{d}(\Sigma_2(\#^nK_*), \s_y)=\sum_{i=1}^n \overline{d}(\Sigma_2(K_*), \s_{y_i}).\]
Therefore, at least one of $\overline{d}(\Sigma_2(\#^nK_*), \s_y)$ and $\overline{d}(\Sigma_2(\#^nK_*), \s_{2y})$ is strictly positive.
\end{proof}

We are now ready to prove our main theorem. We recall the statement.

\begin{reptheorem}{thm:main}There exist topologically slice knots with arbitrarily large concordance $\Z$-genus. More precisely, there exists a topologically slice knot $K_*$ such that for each $n \in \mathbb{N}$ its $n$-fold connected self-sum has concordance $\Z$-genus at least $n$. Furthermore, we can choose  $K_*$ so that it has smooth $4$-ball genus one.
\end{reptheorem}
\begin{proof}Let $K_*$ be the knot in Example~\ref{exl:ourexl}. It is topologically slice and has smooth $4$-ball genus at most one as we observed above. 

We will show that $\#^nK_*$, the $n$-fold connected self-sum of $K_*$, has concordance $\Z$-genus at least $n$. Suppose $\gc(\#^nK_*)=m$ and $n>m$, then by Corollary~\ref{cor:mainobstr} there exists a nontrivial subgroup $N \subseteq H_1(\Sigma_2(\#^nK_*))_5$ such that $\overline{d}(\Sigma_2(\#^nK_*), \s_z) = 0$ for each $z \in N$. This contradicts Proposition~\ref{prop:d-of-2fold}.\end{proof}

\begin{remark}\label{rmk:livingston}
A similar though much more involved argument using Casson-Gordon signatures instead of $d$-invariants can be used to show that the difference between the topological $4$-ball genus and the topological concordance $\Z$-genus can be arbitrarily large.

 In fact, a careful reading of Livingston's paper~\cite{Livingston:2004-1} shows that, despite the fact that it only explicitly considers the topological concordance genus (i.e. the minimal Seifert genus of any representative of a given topological concordance class), all the relevant work has been done
   to prove that for each $n \in \mathbb{N}$ there exists a knot $K_n$ with the  following properties:
\begin{enumerate}
\item $K_n$ has topological $4$-ball genus one. 
\item $K_n$ is not topologically concordant to any knot $J$ with  $r( H_1(\Sigma_2(J))) < 2n$.
\end{enumerate}
In particular, by Proposition~\ref{prop:branchedlower}, this implies that $K_n$ is not topologically concordant to any knot $J$ with  $\gz(J) <n$.
\end{remark}

This result of Livingston is relevant to understanding how $\gt$ changes under the satellite operation. Classical techniques show that for each pattern $P$ there is a constant $g_P\geq g_4(P(U))$ such that  for any knot $K$ we have
\[g_4(P(K)) \leq |w_P| \cdot g_4(K) + g_P,\]
where $w_P$ is the algebraic winding number of the pattern $P$, and modern Heegaard Floer invariants can be used to show that this inequality is sometimes sharp~\cite{Feller-Miller-Pinzon-Caicedo:2019-1}. 

We similarly see that for each $P$ there exists a constant $g_P^{\text{top}}\geq \gt(P(U))$ such that for any knot $K$ we have
\begin{equation}\label{eq:easybound}
\gt(P(K)) \leq |w_P| \cdot \gt(K)+g_P^{\text{top}},\end{equation}
but now the sharpness of this result is far from clear--in fact, it remains possible that the factor coming from winding number is unnecessary. More precisely, it is not known if the following inequality holds for every knot $K$:
\begin{equation}\label{eq:satelliteinequal}\gt(P(K)) \leq \gt(K) + \gt(P(U)). \end{equation}
Work of Feller, Miller, and Pinzon-Caicedo~\cite{Feller-Miller-Pinzon-Caicedo:2019-1} and independently~McCoy~\cite{McCoy:2019-1} shows that for $P$ a pattern with $P(U)=U$ one has
\begin{equation}\label{eq:FMPC} 
\gt(P(K)) \leq \gtc(K),\end{equation}
Hence, Livingston's examples of knots $K$ with $\gt(K) < \gtc(K)$ are good candidates for input companion knots if one wishes to show that the inequality \eqref{eq:satelliteinequal} is not true.  

More concretely, we state the following problem. 

\begin{problem}
Let $K_n$ be Livingston's $n$th knot, which has $\gt(K)=1$ and $\gtc(K) \geq n$. If the inequality~\eqref{eq:satelliteinequal} is true, then for any positive integers  $n, m$ we have that $$\gt(C_{m,1}(K_n))=1.$$ However, the best upper bounds we have coming from the topological concordance $\Z$-genus (i.e.\ the inequality~\eqref{eq:FMPC}) and from classical arguments (i.e.\ the inequality~\eqref{eq:easybound}) are
\[ \gt(C_{m,1}(K_n))\leq  \min\{m, \gtc(K_n)\}.\] 
Determine  $\gt(C_{m,1}(K_n))$ for some $m,n>1$.\end{problem}


\bibliographystyle{alpha}
\def\MR#1{}
\bibliography{bib}
\end{document}